\documentclass{article}

\usepackage{amsmath,amssymb,amsthm}

\usepackage[isolatin]{inputenc}

\newtheorem{theorem}{Theorem}[section]

\theoremstyle{remark}
\newtheorem{remark}{Remark}[section]
\theoremstyle{definition}

\theoremstyle{definition}

\begin{document}

\title{On the Noether Invariance Principle
for Constrained Optimal Control Problems\footnote{Research
report CM04/I-12, Dep. Mathematics, Univ. Aveiro, May 2004.
To be presented at the 6th WSEAS International Conference
on Applied Mathematics, Corfu, Greece, August 17-19, 2004.
Accepted for publication in the journal
\emph{WSEAS Transactions on Mathematics}.}}

\author{Delfim F. M. Torres\\
        \texttt{delfim@mat.ua.pt}\\[0.3cm]
        Department of Mathematics\\
        University of Aveiro\\
        3810-193 Aveiro, Portugal}

\date{\texttt{http://www.mat.ua.pt/delfim}}

\maketitle


\begin{abstract}
We obtain a generalization of Noether's invariance principle
for optimal control problems with equality and inequality
state-input constraints.
The result relates the invariance properties of the problems with
the existence of conserved quantities along the
constrained Pontryagin extremals. A result of this kind
was posed as an open question by Vladimir Tikhomirov, in 1986.
\end{abstract}


\vspace*{0.5cm}



\vspace*{0.3cm}



\section{Introduction}

Noether's invariance principle is one of the most helpful and
fundamental results of physics. It describes the universal fact
that ``invariance with respect to some family of parameter
transformations gives rise to the existence of certain conserved
quantities.'' Such relation is used to explain everything from the
fusion of hydrogen to the motion of planets orbiting the sun
\cite{MR53:10537}. For a modern account of Noether's invariance
principle, in the context of the calculus of variations, we refer
the reader to \cite{MR2000m:49002,brunt}. Extensions for the
unconstrained problems of optimal control are available in
\cite{ejc,torresPortMath}. Here we generalize the previous results
\cite{ejc,torresPortMath} to cover both holonomic and nonholonomic
constraints. The motivation for the present study was Tikhomirov's
book \cite[Sec. 4.3]{MR87m:49004}: \emph{Presumably, it can be
shown for a sufficiently broad class of extremal problems
involving constraints that Noether's invariance theorem is still
valid. But general results of this kind have, to date, not been
obtained.} Theorem~\ref{r:Const:consLaw} provides such general
result.


\section{Constrained Optimal Control Problems and Optimality Conditions}

We deal with a broad class of extremal problems in the calculus
of variations and optimal control involving equality
and/or inequality constraints. We consider a nonlinear control system,
\begin{equation}
\label{eq:ncs}
\dot{x}(t) = \varphi\left(t,x(t),u(t)\right) \, ,
\end{equation}
of $n$ differential equations, subject to $m-m'$ equality constraints,
\begin{equation}
\label{eq:ec:ig}
\phi_i\left(t,x(t),u(t)\right) = 0 \, , \quad i = 1,\ldots, m-m' \, ,
\end{equation}
$m'$ inequality constraints,
\begin{equation}
\label{eq:ec:des}
\phi_j\left(t,x(t),u(t)\right) \ge 0 \, , \quad j = m-m'+1,\ldots, m \, ,
\end{equation}
and $2n$ boundary conditions
\begin{equation}
\label{eq:bc}
x(a) = \alpha \, , \quad x(b) = \beta \, .
\end{equation}
The problem is to find a piecewise-continuous control
function $u(\cdot) = \left(u_1(\cdot),\ldots,u_r(\cdot)\right)$,
and the corresponding state trajectory
$x(\cdot) = \left(x_1(\cdot),\ldots,x_n(\cdot)\right)$, satisfying
\eqref{eq:ncs}, \eqref{eq:ec:ig}, \eqref{eq:ec:des}, and \eqref{eq:bc},
which minimizes or maximizes the integral cost functional
\begin{equation*}
I[x(\cdot),u(\cdot)] = \int_a^b L(t,x(t),u(t)) \mathrm{d}t \, .
\end{equation*}
This problem is denoted in the sequel by $(P)$.
Both the initial time $a$ and terminal time $b$, $a < b$, are fixed.
The boundary values $\alpha,\,\beta \in \mathbb{R}^n$ are also
given. The functions $L(\cdot,\cdot,\cdot)$, $\varphi(\cdot,\cdot,\cdot)$
and $\phi(\cdot,\cdot,\cdot)$ are assumed to be continuously differentiable
with respect to all variables. We shall also assume that the
Jacobian of the constraints \eqref{eq:ec:ig} and \eqref{eq:ec:des},
$\frac{\partial \phi}{\partial u}(t,x,u)$, has
full rank for all $(t,x,u) \in [a,b] \times \mathbb{R}^n \times \mathbb{R}^r$
\cite[Ch. 6]{MR29:3316b}.
Many practical applications of problem $(P)$ appear in
engineering and economics. We refer the interested
reader to \cite{MR93d:49002} and references therein.

The celebrated Pontryagin's maximum principle \cite{MR29:3316b}
gives necessary optimality conditions to be satisfied
by the solutions of optimal control problems.
\begin{theorem}[Pontryagin Maximum Principle for $(P)$]
\label{PMP:P}
Let $u(t)$, $t \in [a,b]$, be an optimal control for
the constrained optimal control problem $(P)$, and
$x(\cdot)$ the corresponding state trajectory. Then there
exists a constant $\psi_0 \le 0$, a continuous costate
$n$-vector function $\psi(\cdot)$ having piecewise-continuous derivatives,
and (assuming that the rank condition is satisfied)
piecewise-continuous multipliers $\lambda(\cdot)$,
$\lambda(t) \ge 0$, satisfying the equations:
\begin{gather*}
\dot{x}(t) =
\frac{\partial H}{\partial \psi}\left(t,x(t),u(t),\psi_0,\psi(t),\lambda(t)\right)
\, , \\
\dot{\psi}(t) = -
\frac{\partial H}{\partial x}\left(t,x(t),u(t),\psi_0,\psi(t),\lambda(t)\right)
\, , \\
\frac{\partial H}{\partial u}\left(t,x(t),u(t),\psi_0,\psi(t),\lambda(t)\right) = 0\, ,
\end{gather*}
along with the constraints \eqref{eq:ec:ig}-\eqref{eq:ec:des}
and boundary conditions \eqref{eq:bc},
where the Hamiltonian $H$ is defined by
\begin{equation*}
H(t,x,u,\psi_0,\psi,\lambda)
= \psi_0 L(t,x,u) + \psi \cdot \varphi(t,x,u) + \lambda \cdot \phi(t,x,u)\, .
\end{equation*}
Moreover, $H(t,x(t),u(t),\psi_0,\psi(t),\lambda(t))$ is a continuous
function of $t$ and, on each interval of continuity of $u(\cdot)$,
is differentiable and satisfies the equality
\begin{equation}
\label{eq:dHdt}
\frac{dH}{dt} = \frac{\partial H}{\partial t} \, .
\end{equation}
\end{theorem}
For versions of Theorem~\ref{PMP:P} under weaker smoothness hypotheses
on the data of the problem, see \cite{MR1897883,MR2001c:49001}.


\section{Main Result}

The following result asserts that
the presence of an invariant structure of the optimal control
problems involving equality and inequality constraints,
imply that their extremals (and solutions) also possess a certain invariance.
The result is expressed, as it happens for the problems
of the calculus of variations \cite{MR2000m:49002,brunt}
and for the unconstrained optimal control problems \cite{ejc,torresPortMath},
as an instance of Noether's universal principle.
Theorem~\ref{r:Const:consLaw} extends
\cite[Theorem 5.1]{torresPortMath} to the case
of constrained optimal control problems.

\begin{theorem}
\label{r:Const:consLaw}
If there exists a $C^2$-smooth one-parameter family of maps
\begin{gather*}
h^s : [a,b] \times \mathbb{R}^n \times \mathbb{R}^r \rightarrow
       \mathbb{R} \times \mathbb{R}^n \times \mathbb{R}^r \, , \\
h^s(t,x,u) = \left(T(t,x,u,s), X(t,x,u,s), U(t,x,u,s)\right) \, , \\
s \in (-\varepsilon, \varepsilon)  \, , \, \varepsilon > 0 \, ,
\end{gather*}
with $h^0(t,x,u) = (t,x,u)$ for all
$(t,x,u) \in [a,b] \times \mathbb{R}^n \times \mathbb{R}^r$, and
satisfying
\begin{gather}
L\left(t,x(t),u(t)\right)
= L \circ h^s\left(t,x(t),u(t)\right)
\frac{d}{dt} T\left(t,x(t),u(t),s\right) \, , \label{eq:invi} \\
\frac{d}{dt} X\left(t,x(t),u(t),s\right)
= \varphi \circ h^s\left(t,x(t),u(t)\right)
\frac{d}{dt} T\left(t,x(t),u(t),s\right) \, ,  \label{eq:invii} \\
\phi\left(t,x(t),u(t)\right)
= \phi \circ h^s\left(t,x(t),u(t)\right)
\frac{d}{dt} T\left(t,x(t),u(t),s\right)  \, , \label{eq:inviii}
\end{gather}
then,
\begin{equation*}
\psi(t) \cdot \frac{\partial}{\partial s}
\left.X\left(t,x(t),u(t),s\right)\right|_{s = 0}
- H(t,x(t),u(t),\psi_0,\psi(t),\lambda(t)) \frac{\partial}{\partial s}
\left.T\left(t,x(t),u(t),s\right)\right|_{s = 0}
\end{equation*}
is constant in $t \in [a,b]$
for any quintuple $\left(x(\cdot),u(\cdot),\psi_0,\psi(\cdot),\lambda(\cdot)\right)$
satisfying the Pontryagin maximum principle (Theorem~\ref{PMP:P}),
with $H$ the Hamiltonian associated
to the problem $(P)$: $H(t,x,u,\psi_0,\psi,\lambda)
= \psi_0 L(t,x,u) + \psi \cdot \varphi(t,x,u) + \lambda \cdot \phi(t,x,u)$.
\end{theorem}

\begin{proof}
Using the fact that $h^0(t,x,u) = (t,x,u)$,
from condition \eqref{eq:invi} one gets
\begin{align}
0 & = \frac{d}{ds}
\left.\left(L \circ h^s\left(t,x(t),u(t)\right)
\frac{d}{dt} T\left(t,x(t),u(t),s\right)\right)\right|_{s = 0} \notag \\
  & = \frac{\partial L}{\partial t}
      \left.\frac{\partial T}{\partial s}\right|_{s = 0}
      + \frac{\partial L}{\partial x} \cdot
      \left.\frac{\partial X}{\partial s}\right|_{s = 0}
      + \frac{\partial L}{\partial u} \cdot
      \left.\frac{\partial U}{\partial s}\right|_{s = 0}
      + L \frac{d}{dt} \left.\frac{\partial T}{\partial s}\right|_{s =
      0} \, , \label{eq:ds0i}
\end{align}
while condition \eqref{eq:invii} and \eqref{eq:inviii} yields
\begin{gather}
\frac{d}{dt} \left.\frac{\partial X}{\partial s}\right|_{s = 0}
= \frac{\partial \varphi}{\partial t}
      \left.\frac{\partial T}{\partial s}\right|_{s = 0}
      + \frac{\partial \varphi}{\partial x} \cdot
      \left.\frac{\partial X}{\partial s}\right|_{s = 0}
      + \frac{\partial \varphi}{\partial u} \cdot
      \left.\frac{\partial U}{\partial s}\right|_{s = 0}
      + \varphi \frac{d}{dt} \left.\frac{\partial T}{\partial s}\right|_{s =
      0} \, , \label{eq:ds0ii} \\
0 = \frac{\partial \phi}{\partial t}
      \left.\frac{\partial T}{\partial s}\right|_{s = 0}
      + \frac{\partial \phi}{\partial x} \cdot
      \left.\frac{\partial X}{\partial s}\right|_{s = 0}
      + \frac{\partial \phi}{\partial u} \cdot
      \left.\frac{\partial U}{\partial s}\right|_{s = 0}
      + \phi \frac{d}{dt} \left.\frac{\partial T}{\partial s}\right|_{s =
      0} \, . \label{eq:ds0iii}
\end{gather}
Multiplying \eqref{eq:ds0i} by $\psi_0$, \eqref{eq:ds0ii} by $\psi(t)$,
and \eqref{eq:ds0iii} by $\lambda(t)$, we can write:
\begin{multline}
\label{eq:joined}
\psi_0 \left(\frac{\partial L}{\partial t}
      \left.\frac{\partial T}{\partial s}\right|_{s = 0}
      + \frac{\partial L}{\partial x} \cdot
      \left.\frac{\partial X}{\partial s}\right|_{s = 0}
      + \frac{\partial L}{\partial u} \cdot
      \left.\frac{\partial U}{\partial s}\right|_{s = 0}
      + L \frac{d}{dt} \left.\frac{\partial T}{\partial s}\right|_{s =
      0}\right) \\
+ \psi(t) \cdot \left(\frac{\partial \varphi}{\partial t}
      \left.\frac{\partial T}{\partial s}\right|_{s = 0}
      + \frac{\partial \varphi}{\partial x} \cdot
      \left.\frac{\partial X}{\partial s}\right|_{s = 0}
      + \frac{\partial \varphi}{\partial u} \cdot
      \left.\frac{\partial U}{\partial s}\right|_{s = 0}
      + \varphi \frac{d}{dt} \left.\frac{\partial T}{\partial s}\right|_{s =
      0} - \frac{d}{dt} \left.\frac{\partial X}{\partial s}\right|_{s = 0}\right)\\
+ \lambda(t) \cdot \left(\frac{\partial \phi}{\partial t}
      \left.\frac{\partial T}{\partial s}\right|_{s = 0}
      + \frac{\partial \phi}{\partial x} \cdot
      \left.\frac{\partial X}{\partial s}\right|_{s = 0}
      + \frac{\partial \phi}{\partial u} \cdot
      \left.\frac{\partial U}{\partial s}\right|_{s = 0}
      + \phi \frac{d}{dt} \left.\frac{\partial T}{\partial s}\right|_{s =
      0} \right) = 0 \, .
\end{multline}
According to the Pontryagin maximum principle, the function
\begin{multline*}
\psi_0 L\left(t,x(t),U\left(t,x(t),u(t),s\right)\right)
+ \psi(t) \cdot \varphi\left(t,x(t),U\left(t,x(t),u(t),s\right)\right)\\
+ \lambda(t) \cdot \phi\left(t,x(t),U\left(t,x(t),u(t),s\right)\right)
\end{multline*}
attains an extremum for $s = 0$. Therefore
\begin{equation*}
\psi_0 \frac{\partial L}{\partial u} \cdot
\left.\frac{\partial U}{\partial s}\right|_{s = 0}
+ \psi(t) \cdot \frac{\partial \varphi}{\partial u} \cdot
\left.\frac{\partial U}{\partial s}\right|_{s = 0}
+ \lambda(t) \cdot \frac{\partial \phi}{\partial u} \cdot
\left.\frac{\partial U}{\partial s}\right|_{s = 0} = 0
\end{equation*}
and \eqref{eq:joined} simplifies to
\begin{multline}
\label{eq:applASandDHdt}
\psi_0 \left(\frac{\partial L}{\partial t}
      \left.\frac{\partial T}{\partial s}\right|_{s = 0}
      + \frac{\partial L}{\partial x} \cdot
      \left.\frac{\partial X}{\partial s}\right|_{s = 0}
      + L \frac{d}{dt} \left.\frac{\partial T}{\partial s}\right|_{s =
      0}\right) \\
+ \psi(t) \cdot \left(\frac{\partial \varphi}{\partial t}
      \left.\frac{\partial T}{\partial s}\right|_{s = 0}
      + \frac{\partial \varphi}{\partial x} \cdot
      \left.\frac{\partial X}{\partial s}\right|_{s = 0}
      + \varphi \frac{d}{dt} \left.\frac{\partial T}{\partial s}\right|_{s =
      0} - \frac{d}{dt} \left.\frac{\partial X}{\partial s}\right|_{s = 0}\right)\\
+ \lambda(t) \cdot \left(\frac{\partial \phi}{\partial t}
      \left.\frac{\partial T}{\partial s}\right|_{s = 0}
      + \frac{\partial \phi}{\partial x} \cdot
      \left.\frac{\partial X}{\partial s}\right|_{s = 0}
      + \phi \frac{d}{dt} \left.\frac{\partial T}{\partial s}\right|_{s =
      0} \right) = 0 \, .
\end{multline}
From the adjoint system
$\dot{\psi} = -\frac{\partial H}{\partial x}$ and the equality
\eqref{eq:dHdt}, we know that
\begin{gather*}
\dot{\psi} = - \psi_0 \frac{\partial L}{\partial x}
- \psi \cdot \frac{\partial \varphi}{\partial x}
- \lambda \cdot \frac{\partial \phi}{\partial x} \\
\frac{d}{dt} H = \psi_0 \frac{\partial L}{\partial t}
+ \psi \cdot \frac{\partial \varphi}{\partial t}
+ \lambda \cdot \frac{\partial \phi}{\partial t} \, ,
\end{gather*}
and one concludes that \eqref{eq:applASandDHdt} is equivalent to
\begin{equation*}
\frac{d}{dt} \left(\psi(t) \cdot
\left.\frac{\partial X}{\partial s}\right|_{s = 0}
- H \left.\frac{\partial T}{\partial s}\right|_{s = 0}\right)
= 0 \, .
\end{equation*}
The proof is complete.
\end{proof}

\begin{remark}
Theorem~\ref{r:Const:consLaw} is still valid in the situation
where the boundary values of the state variables are not fixed.
We have considered conditions \eqref{eq:bc} only to simplify
the presentation of the maximum principle: transversality
conditions are not relevant in the proof of our result.
\end{remark}
\begin{remark}
It is possible to give a formulation of
Theorem~\ref{r:Const:consLaw} under gauge-variance
(see \cite{MR83c:70020} for the
concept of gauge-variance), and deal with equalities
\eqref{eq:invi}--\eqref{eq:inviii} up to first-order
terms in the parameter $s$ (see the quasi-invariance notion
introduced in \cite{torresPortMath} for the unconstrained
optimal control problem).
\end{remark}

We now illustrate the application of Theorem~\ref{r:Const:consLaw}
with an example.


\section{Extraction of an Exhaustible Resource}

We borrow from \cite[pp.~194--198]{MR93d:49002} a simple
example of an optimal control problem in economics with
one state variable, two control variables, and an equality
constraint on the state and control variables
($n=1$, $r=2$, $m=1$, $m'=0$):
\begin{gather}
\int_0^T u_1^\gamma(t)  \mathrm{d}t \rightarrow \max \notag \\
\dot{x}(t) = -u_2(t) \, , \label{eq:prb:eer}\\
x^{\alpha \gamma}(t) u_2^{\beta \gamma}(t) - u_1^\gamma(t) = 0 \, , \notag \\
x(0) = x_0 \, , \quad x(T) = x_T \, , \notag
\end{gather}
where $\gamma < 1$, $\alpha + \beta < 1$ ($\alpha, \beta > 0$),
and $x_T < x_0$. Here $x(t)$ denote the stock of an exhaustible
resource at time $t$; $u_2(\cdot)$ the rate of extraction from the
stock; and $u_1(\cdot)$ the flow of consumption of the finished good:
see \cite{MR93d:49002} for more details on the model, and for an
economic interpretation of the Pontryagin maximum principle.
Problem \eqref{eq:prb:eer} is invariant under the one-parameter
family of transformations $h^s(t,x,u_1,u_2)
= \left(T(t,s), X(x,s), U_1(u_1,s), U_2(u_2,s)\right)$ defined by
\begin{gather}
T(t,s) = \mathrm{e}^{-\gamma(\alpha+\beta)s}t \, , \quad
X(x,s) = \mathrm{e}^{(1-\beta\gamma)s}x \, , \label{eq:trfEq:ex} \\
U_1(u_1,s) = \mathrm{e}^{(\alpha+\beta)s}u_1 \, , \quad
U_2(u_2,s) = \mathrm{e}^{(\alpha\gamma+1)s}u_2 \, , \notag
\end{gather}
which coincides, for $s = 0$, with the identity transformation:
$h^0(t,x,u_1,u_2) = \left(t,x,u_1,u_2\right)$. In fact, for problem
\eqref{eq:prb:eer} one has $L(u_1) = u_1^\gamma$,
$\varphi(u_2) = -u_2$, and
$\phi(x,u_1,u_2) = x^{\alpha\gamma} u_2^{\beta\gamma} - u_1^\gamma$,
and all conditions \eqref{eq:invi}, \eqref{eq:invii}, and \eqref{eq:inviii}
are satisfied under \eqref{eq:trfEq:ex}:
\begin{gather*}
\begin{split}
L(U_1) \frac{d}{dt} T(t,s) &= \mathrm{e}^{\gamma(\alpha+\beta)s}u_1^\gamma
\mathrm{e}^{-\gamma(\alpha+\beta)s} = u_1^\gamma \\
&= L(u_1) \, ,
\end{split} \\
\begin{split}
\varphi(U_2) \frac{d}{dt} T(t,s) &= - \mathrm{e}^{(\alpha\gamma+1)s}u_2
\mathrm{e}^{-\gamma(\alpha+\beta)s} =
- \mathrm{e}^{(1-\beta\gamma)s} u_2 =
\mathrm{e}^{(1-\beta\gamma)s} \dot{x} \\
&= \frac{d}{dt} X(x,s) \, ,
\end{split} \\
\begin{split}
\phi(X,U_1,U_2) \frac{d}{dt} T(t,s) &=
\left(\mathrm{e}^{\alpha\gamma(1-\beta\gamma)s}x^{\alpha\gamma}
\mathrm{e}^{\beta\gamma(\alpha\gamma+1)s}u_2^{\beta\gamma}
- \mathrm{e}^{\gamma(\alpha+\beta)s}u_1^\gamma\right)
\mathrm{e}^{-\gamma(\alpha+\beta)s}\\
&= \mathrm{e}^{\gamma(\alpha+\beta)s}
\left(x^{\alpha\gamma} u_2^{\beta\gamma} - u_1^\gamma\right)
\mathrm{e}^{-\gamma(\alpha+\beta)s}\\
&= \phi(x,u_1,u_2) \, .
\end{split}
\end{gather*}
Having in mind that the problem is autonomous, that is,
the Hamiltonian $H(t,x,u_1,u_2,\psi_0,\psi,\lambda)
= \psi_0 u_1^\gamma - \psi u_2
+ \lambda \left(x^{\alpha\gamma} u_2^{\beta\gamma} - u_1^\gamma\right)$
does not depend on $t$, and that from equality \eqref{eq:dHdt}
this implies the Hamiltonian $H$ to be constant along the extremals,
it follows from our Theorem~\ref{r:Const:consLaw} that
\begin{equation}
\label{eq:leiCons:Ex}
(1-\beta\gamma) \psi(t) x(t) + \gamma H (\alpha+\beta) t = constant \, .
\end{equation}
Relation \eqref{eq:leiCons:Ex} gives an immediate insight
about the solution of the problem -- the rate of the
value times the stock size of the resource must
be constant: $\frac{d}{dt} \left(\psi(t) x(t)\right) = const$.
This conservation law has an economical interpretation
in terms of the Cobb-Douglas form of the production function.


\section*{Acknowledgements}

Research supported by the Portuguese Foundation for Science and Technology (FCT),
partially via the FEDER Project POCTI/MAT/41683/2001, partially via the R\&D unit
Centre for Research in Optimization and Control (CEOC) of the University of Aveiro.



\end{document}